\documentclass[12pt]{article}
\pdfoutput=1
\usepackage[margin=30pt,font=small,labelfont=bf, labelsep=endash]{caption}
\usepackage{authblk}
\usepackage{amsmath}
\usepackage{amssymb}
\usepackage{url}
\usepackage{amsthm}
\usepackage{mathtools}
\DeclarePairedDelimiter{\ceil}{\lceil}{\rceil}
\usepackage{graphicx}
\usepackage{fullpage}
\usepackage[parfill]{parskip}

\theoremstyle{definition}
\newtheorem{definition}{Definition}
\newtheorem{theorem}{Theorem}

\newtheorem{observation}{Observation}
\newtheorem{conjecture}{Conjecture}

\begin{document}

\title{Strengthening strong immersions with Kempe chains}
\author[1]{Todd A. Gibson}
\author{
Todd A. Gibson\thanks{Corresponding Author. Email: tagibson@csuchico.edu}
\\
Department of Computer Science\\
California State University, Chico\\
Chico, CA 95929, USA\\
}
\maketitle

\begin{abstract}
  Every properly colored graph with $\chi(G)=k$ colors
  has edge-disjoint Kempe \mbox{``backbones''}, Kempe chains anchored by color-critical vertices for each pair of colors.  Certain color permutations arrange these backbones into a clique-like structure, a strengthening of strong immersions of complete graphs. This strengthened immersion is suggested as a template for identifying the disjoint subgraphs comprising Hadwiger's conjectured $K_k$ minor present in $k$-chromatic graphs.
\end{abstract}

\section{Introduction}
An immersion of the complete graph $K_k$ in a simple, undirected graph $G$ is an injective function that maps $f:V(K_k)\rightarrow{V(G)}$. Each edge $(u,v)\in E(K_k)$ corresponds to a path in $G$ with endpoints $(f(u),f(v))$. The immersion is strong if the paths are internally disjoint. Studying strong immersions of complete graphs has been motivated by their potential to make progress towards resolving Hadwiger's conjecture~\cite{Abu-khzam2003}. Here we identify additional structural elements of graphs closely aligned with Kempe chains and how they assemble into a strengthened form of $K_k$ immersions. These \emph{Kempe cliques} are attractive research targets because of their close alignment to the graph's chromatic number.

Let $G$ be a simple, undirected graph with chromatic number $\chi(G)=k$, and chromatic coloring, that is, properly colored with $C={c_1, c_2, \dots, c_k}$ colors.

\begin{figure}[h]
  \begin{center}\scalebox{1.30}{\includegraphics{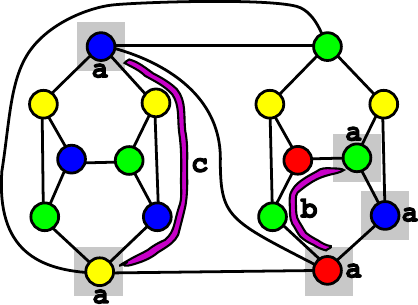}}\end{center}
    \caption{Annotated structures. Critical vertices, \textbf{a}, each framed by a grey square. A (green,~red)-Kempe backbone, \textbf{b}, alongside a purple ribbon. A (blue, yellow)-Kempe backbone, \textbf{c}, alongside a purple ribbon. A length-3 (yellow, green)-Kempe backbone is present but not highlighted. Three length-1 Kempe backbones, (red, yellow), (red, blue), and (blue, green) complete the list of six Kempe backbones present in a $\chi(G)=4$ graph.  \label{fig:graph_definitions}} \end{figure}

\definition A \textit{critical vertex}, $X_i$, is a vertex with color $c_i$ that is adjacent to at least $k-1$ neighbors, colored with all remaining colors in $C$ (Figure~\ref{fig:graph_definitions}).

\definition A \textit{Kempe chain} contains vertex $v$ colored $c_i$, and is the maximal connected subgraph of vertices colored either $c_i$ or $c_j$. Also referred to as a $(c_i, c_j)$-Kempe chain.

\definition A \textit{Kempe backbone} is a path within a ($c_i$, $c_j$)-Kempe chain leading from critical vertex $X_i$ to critical vertex $X_j$. These critical vertices are the \textit{anchors} 
of the backbone. 
Also referred to as a ($c_i, c_j)$-Kempe backbone (Figure~\ref{fig:graph_definitions}).

\definition A \textit{Kempe swap} interchanges the colors of a $(c_i, c_j)$-Kempe chain. Each vertex in the chain assumes the color of its chain neighbors.

\section{Kempe backbones}
\observation \label{critverts} Graph $G$ contains at least $k$ critical vertices. Each color in $C$ labels at least one critical vertex.
\proof Assume that $G$ is missing critical vertex $X_a$. Each $c_a$ colored vertex may be assigned a different color because no $c_a$ colored vertex is adjacent to all remaining colors by the definition of a critical vertex. Therefore, $G$ is $(k-1)$-colorable, a contradiction.\qed

\observation \label{swappreserve} A Kempe swap preserves the criticality (whether a vertex is critical or non-critical) of each vertex in the Kempe chain.
\proof Each vertex in the $(c_a, c_b)$-Kempe chain is adjacent to neighbors labeled with $1\le n<k$ colors. Interchanging the color of vertex $v$ in the chain from $c_a$ to $c_b$ replaces $c_b$-labeled neighbors with $c_a$-labeled neighbors (and visa versa). Other neighbors are not affected by the swap.  Therefore, $n$ does not change for each vertex participating in the swap. \qed

Note that although the criticality is preserved for chain members undergoing a swap, the criticality of vertices adjacent to members of the swapped chain may change.

The following theorem describes the Kempe backbone characteristic of $G$. There are $(k^2-k)/2$ unique pairs of colors. For each unique color pair $(c_i,c_j)$, there exist two critical vertices, $(X_i, X_j)$, anchoring a $(c_i,c_j)$-Kempe backbone.

\theorem \label{allpairs} For each pair of colors $(c_i,c_j)$ with $c_i, c_j \in C$, $i \ne{j}$
there is a critical vertex pair $(X_i, X_j)$ that anchors a $(c_i, c_j)$-Kempe backbone.
\proof Let $(a,b)$ be such a color pair and assume instead that no $(X_a, X_b)$ pair of critical vertices anchors an $(a,b)$-Kempe backbone. By Observation~\ref{critverts}, critical vertices $X_a$ and $X_b$ exist. 
So, $X_a$ and $X_b$ are members of different, disconnected $(a,b)$-Kempe chains. Swap the colors of each $(a,b)$-Kempe chain containing a critical vertex $X_a$. By Observation~\ref{swappreserve}, this replaces each $X_a$ with a critical vertex of color $b$. This process may be repeated to eliminate all critical vertices of color $a$, which is a contradiction by Observation~\ref{critverts}.\qed

Graph $G$ may have more than one critical vertex of a given color; not all critical vertices are required to participate in a Kempe backbone to satisfy Theorem~\ref{allpairs}. Also, more than one $(c_a,c_b)$-Kempe backbone may exist in $G$.

\begin{figure}[t]
  \begin{center}\scalebox{1.30}{\includegraphics{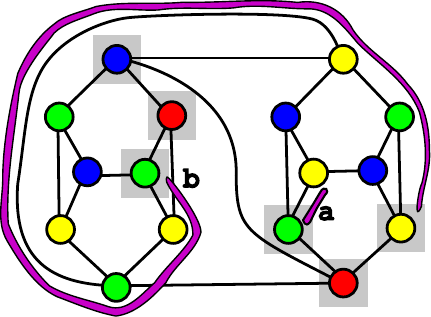}}\end{center}
    \caption{A chromatically colored $k=4$ graph. Critical vertices are framed
    by grey squares.  Swapping the (green, yellow)-Kempe chain, \textbf{a}, 
    will eliminate the green critical vertex. The second (green,
    yellow)-Kempe chain, \textbf{b}, includes a (green, yellow)-Kempe backbone
    following the purple ribbon. Swapping this chain merely interchanges the critical
    vertices of the backbone. \label{fig:illustrate_thm1}} 
  \end{figure}

Theorem~\ref{allpairs} describes the architecture that prevents G from using fewer than $k$ colors. Any attempt to use $(a,b)$-Kempe swaps to remove all $X_a$ critical vertices will encounter at least one $X_a$ that has a Kempe backbone to $X_b$ (Figure~\ref{fig:illustrate_thm1}). A $(c_a, c_b)$-Kempe swap on the backbone anchored by $X_a$ and $X_b$ will merely interchange the critical vertices, failing to eliminate $X_a$.

\newpage
\section{Kempe cliques}
\definition A \textit{correctly colored} graph is properly colored with $q\ge\chi(G)$ colors such that it includes $q$ critical vertices, $\{X_1, X_2, \dots, X_q\}$, each labeled with a different color, anchoring $(q^2-q)/2$ Kempe backbones. 
  \definition A \textit{Kempe clique}, $Q_q$, is the collection of $(q^2-q)/2$ Kempe backbones anchored by $q$ critical vertices in a correct coloring.

\begin{figure}[h]
  \begin{center}\scalebox{1.30}{\includegraphics{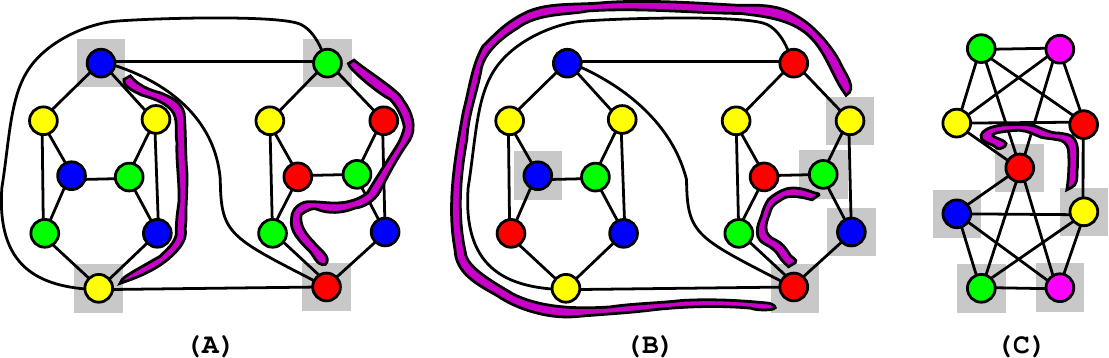}}\end{center}
    \caption{Examples of correct coloring. Critical vertices are framed by grey squares. Kempe backbones greater than length-1 follow purple ribbons. (A) A graph with $\chi(G)=4$.  The Kempe clique includes four length-1 Kempe backbones, one of length-3 and one of length-5. (B) A graph with $\chi(G)=4$. In addition to the Kempe backbones comprising a Kempe clique, a correctly colored graph may have additional critical vertices and additional Kempe backbones. (C) A correctly colored graph with $\chi(G)=5$. The Kempe clique includes a length-3 Kempe backbone.\label{fig:correctly_colored}} \end{figure}

The presence of a Kempe clique does not directly follow from Theorem~\ref{allpairs}. Figures~\ref{fig:graph_definitions}~and~\ref{fig:illustrate_thm1} are examples of proper (and chromatic) colorings that are not correct colorings.  A different color permutation of a properly colored graph may be required to reveal a Kempe clique.  Figure~\ref{fig:correctly_colored} shows labeled examples of correctly colored graphs and their Kempe cliques.

\begin{figure}[t]
  \begin{center}\scalebox{1.30}{\includegraphics{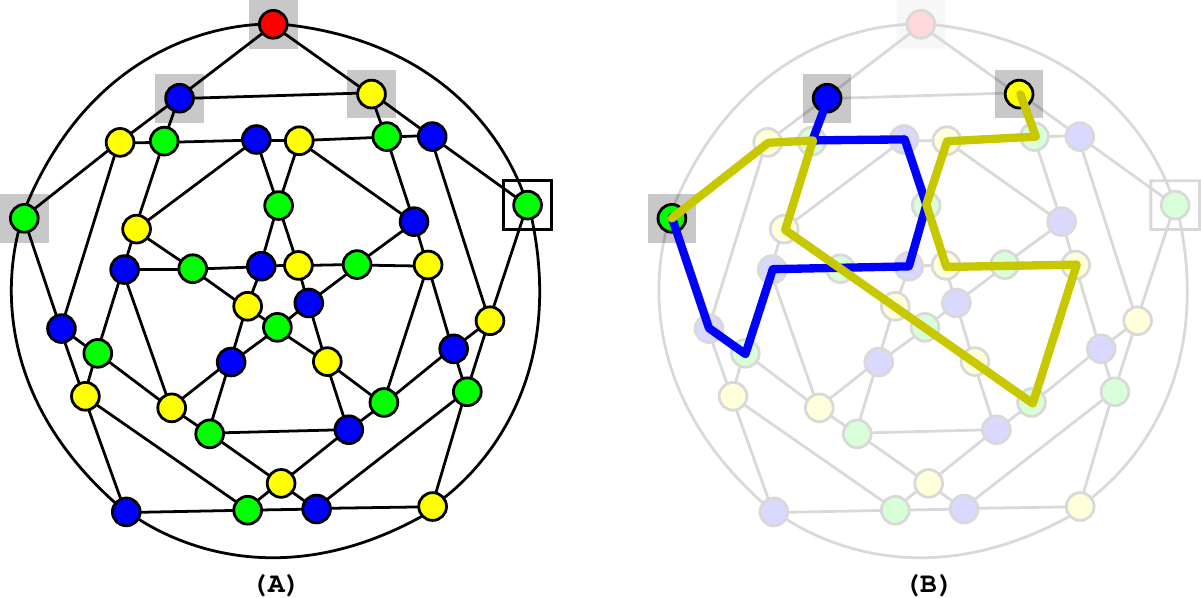}}\end{center}
    \caption{(A) A correct coloring of Koester's graph. (B) Highlighted are the Kempe backbones with length $l>1$.\label{fig:koester_graph}} \end{figure}

 Kempe cliques strengthen strong immersions of complete graphs~\cite{Lescure1988,Devos2014}. Every Kempe clique, $Q_q$, is a strong immersion of the complete graph, $K_q$, but the converse does not follow.  
 For example, Koester's planar graph~\cite{Koester1985}  with $\chi(G)=4$ (Figure~\ref{fig:koester_graph}) has a strongly immersed $K_5$, as do all 4-regular graphs~\cite{Lescure1988}. However, there is no correct coloring of Koester's graph when properly colored with $q=5$ colors. It is straightforward to show that $Q_5$ is non-planar. In fact, a proof that all $\chi(G)=5$ graphs admit a correct coloring would serve as an alternate proof of the four-color theorem. 

\conjecture \label{conjecture:kempeclique} Every simple, undirected graph can be correctly colored.

\begin{figure}[t]
  \begin{center}\scalebox{1.30}{\includegraphics{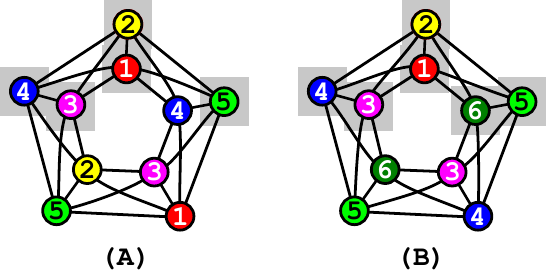}}\end{center}
    \caption{Catlin's counterexample construction, the crossproduct of $C_5$ and $K_2$, $\chi(G)=5$. (A) Graph $G$ is chromatically-colored. Although all vertices are critical, only a subset are framed in grey to ease identification of Kempe backbones. No correct coloring is possible with 5 colors; there is no set of 10 Kempe backbones anchored by 5 critical vertices. (B) Graph $G$ is correctly colored with $6$ colors.\label{fig:catlin2_2}} 
  \end{figure}

Note that a correct coloring uses $q\ge\chi(G)$ colors; $G$ need only be properly-colored, not chromatically-colored. An example of a correct coloring requiring $q>\chi(G)$ can be found in Catlin, 1979~\cite{Catlin1979}. Catlin constructs a counterexample to the Haj\'os Conjecture that every simple graph has a $K_k$ subdivision by taking the crossproduct of a cycle and complete graph (Figure~\ref{fig:catlin2_2}). Given cycle length $2n+1$ and complete graph order $k$, the construction's chromatic number is $\chi(G)=2k+\ceil{k/n}$~\cite{Catlin1979}. A correct coloring of the same construction requires $q=3k$ colors.

A proof of Conjecture~\ref{conjecture:kempeclique} would allow the possibility that the Kempe clique itself serves as the template for the $K_k$ minor conjectured by Hadwiger~\cite{Hadwiger1943}. In particular, each critical vertex of the Kempe clique ``seeds'' a different, disjoint connected subgraph that identifies to form a complete graph. For example, the critical vertices in Figure~\ref{fig:catlin2_2}B each seed a different, disjoint subgraph that forms a $K_6$ minor. Critical vertices $1$ and $2$ are sole members of their subgraphs. Critical vertices $3$, $4$, $5$, and $6$ each include a second vertex in their respective subgraphs.

\conjecture \label{cliqueminor} 
The critical vertices of the Kempe clique in a correctly colored graph with $q$ colors are the seeds of a $K_q$ minor. That is, each is a member of a different, disjoint connected subgraph. Identifying each subgraph forms $K_q$.  Given $h$, the Hadwiger number, $\chi(G) \le q \le h$. 

There are some graph families that are known to be correctly colorable.

\theorem\label{critKminus1}
A $k$-critical graph with minimum degree $\delta(G)=k-1$ admits a correct coloring.
\proof   Select vertex $v$ with $k-1$ neighbors. Since $G-v$ is $(k-1)$-colorable, 
$v$ is the only vertex labeled $c_a$. By Observation~\ref{critverts}, $v$ is a critical vertex. 
By the same Observation, the graph must contain $k-1$ additional critical vertices labeled with the remaining $C-c_a$ colors. To be critical, each of these $k-1$ additional critical vertices must be adjacent to a $c_a$-labeled vertex. Therefore, the $k-1$ neighbors of $v$ comprise the remaining critical vertices in the graph. Because the graph has only $k$ critical vertices, these comprise all anchors for the graph's Kempe backbones, forming a Kempe clique.\qed 

In Abu-khzam and Langston, 2003~\cite{Abu-khzam2003}, Corollary 2 similarly identifies an immersed $K_k$ in a color critical graph with $\delta(G)=k-1$.

    \theorem Uniquely-colorable graphs are correctly colored. 
    \proof Because the subgraph induced by the union of two color classes is connected in a uniquely-colored graph~\cite{harary1969uniquely}, the Kempe backbones form a Kempe clique.  \qed
    
\bibliography{4color}

\begin{thebibliography}{1}

\bibitem{Abu-khzam2003}
Faisal~N Abu-khzam and Michael~A Langston.
\newblock {Graph coloring and the immersion order}.
\newblock {\em Lecture Notes in Computer Science}, 2697:394--403, 2003.

\bibitem{Catlin1979}
Paul~A Catlin.
\newblock {Haj{\'{o}}s' graph-coloring conjecture: Variations and
  counterexamples}.
\newblock {\em Journal of Combinatorial Theory, Series B}, 26(2):268--274, apr
  1979.

\bibitem{Devos2014}
Matt Devos, Zden{\v{e}}k Dvo$\backslash$vr{\'{a}}k, Jacob Fox, Jessica
  McDonald, Bojan Mohar, and Diego Scheide.
\newblock {A minimum degree condition forcing complete graph immersion}.
\newblock {\em Combinatorica}, 34(3):279--298, jun 2014.

\bibitem{Hadwiger1943}
Hugo Hadwiger.
\newblock {{\"{U}}ber Eine Klassifikation Der Streckenkomplexe}.
\newblock {\em Vierteljschr. Naturforsch. Ges. Z{\"{u}}rich}, 88:133--142,
  1943.

\bibitem{harary1969uniquely}
Frank Harary, S~T Hedetniemi, and R~W Robinson.
\newblock {Uniquely colorable graphs}.
\newblock {\em Journal of Combinatorial Theory}, 6(3):264--270, 1969.

\bibitem{Koester1985}
G~Koester.
\newblock {Note to a problem of T. Gallai and G. A. Dirac}.
\newblock {\em Combinatorica}, 5(3):227--228, 1985.

\bibitem{Lescure1988}
F.~Lescure and H.~Meyniel.
\newblock {On a Problem upon Configurations Contained in Graphs with Given
  Chromatic Number}.
\newblock {\em Annals of Discrete Mathematics}, 41(C):325--331, 1988.

\end{thebibliography}

\bibliographystyle{plain}

\end{document}